\documentclass[12pt,twoside]{article}
\usepackage{amsmath}
\usepackage{amsfonts}
\usepackage{amsthm}
\usepackage{amssymb}
\usepackage{graphicx}
\usepackage{multicol}

\begin{document}
\title{{\normalsize
{\bf Ribbonness of a stable-ribbon surface-link, II. General case}}}
\author{{\footnotesize Akio Kawauchi}\\
\date{}
{\footnotesize{\it Osaka Central Advanced Mathematical Institute, Osaka Metropolitan University}}\\ 
{\footnotesize{\it Sugimoto, Sumiyoshi-ku, Osaka 558-8585, Japan}}\\ 
{\footnotesize{\it kawauchi@omu.ac.jp}}}
\date\, 
\date\, 
\maketitle
\vspace{0.25in}
\baselineskip=10pt
\newtheorem{Theorem}{Theorem}[section]
\newtheorem{Conjecture}[Theorem]{Conjecture}
\newtheorem{Lemma}[Theorem]{Lemma}
\newtheorem{Sublemma}[Theorem]{Sublemma}
\newtheorem{Proposition}[Theorem]{Proposition}
\newtheorem{Corollary}[Theorem]{Corollary}
\newtheorem{Claim}[Theorem]{Claim}
\newtheorem{Definition}[Theorem]{Definition}
\newtheorem{Example}[Theorem]{Example}

\begin{abstract} It is shown that any handle-irreducible summand of every stable-ribbon surface-link is a unique ribbon surface-link up to equivalences, so that every stable-ribbon surface-link is a ribbon surface-link. This is a generalization of a previously observed result for a stably trivial surface-link. Two observations are given. One  observation is that a connected sum of two surface-links is a ribbon surface-link if and only if both the connected summands are ribbon surface-links. The other observation is a characterization of when a surface-link consisting of ribbon surface-knot components becomes a ribbon surface-link.

\phantom{x}

\noindent{\it Keywords}: Ribbon surface-link,\, Stable-ribbon surface-link,\, Trivial surface-link,\, Connected sum,\,
O2-handle pair,\, Uniqueness, SUPH system.

\noindent{\it 2020 Mathematics Subject Classification}: Primary 57K45; Secondary 57K40

\end{abstract}

\baselineskip=15pt

\bigskip

\noindent{\bf 1. Introduction}

This paper generalizes the previous result that a stably trivial surface-link is a trivial surface-link to the result that a stable-ribbon surface-link is a ribbon surface-link, \cite{K21}, \cite{K23}.
A {\it surface-link} is a closed oriented 
(possibly disconnected) surface $F$ which is embedded in the 4-space ${\mathbf R}^4$ 
by a smooth  embedding. 
When $F$ is connected, it is also called a {\it surface-knot}. 
When a fixed (possibly disconnected) closed surface $\mathbf F$ is smoothly embedded 
into ${\mathbf R}^4$, it is also 
called an $\mathbf F$-{\it link}. If $\mathbf F$ is the disjoint union of 
some copies of the 2-sphere $S^2$, then it is also called an $S^2$-{\it link}. 
When $\mathbf F$ is connected, it is also called an $\mathbf F$-{\it knot}, 
and an $S^2$-{\it knot} for ${\mathbf F}=S^2$. 
Two surface-links $F$ and $F'$ are {\it equivalent} by an {\it equivalence} $f$ 
if $f$ is an orientation-preserving diffeomorphism  $f:{\mathbf R}^4\to {\mathbf R}^4$ 
sending  $F$  to $F'$ that preserves orientation. 
A {\it trivial } surface-link is a surface-link $F$ which bounds disjoint 
handlebodies smoothly embedded in ${\mathbf R}^4$, where a handlebody is a 3-manifold which is a 3-ball, solid torus or a disk sum of some number of solid tori. 
A trivial surface-knot is also called an {\it unknotted } surface-knot. 
A trivial disconnected surface-link is also called an {\it unknotted-unlinked} 
surface-link. For any given closed oriented 
(possibly disconnected) surface ${\mathbf F}$, a trivial ${\mathbf F}$-link exists uniquely up to equivalences (see \cite{HoK}). 
A {\it ribbon} surface-link is a surface-link $F$ which is obtained 
from a trivial $n S^2$-link $O$ for some $n$ (where $n S^2$ denotes 
the disjoint union of $n$ copies of the 2-sphere $S^2$) 
by surgery along an embedded 1-handle system,  \cite{KSS-2}, \cite{Yajima1}, 
\cite{Yajima2}, \cite{Yana}.    
This object is an old concept in surface-knot theory, but in recent years it is considered as a chord diagram which is a relaxed version of a virtual graph  including a virtual knotoid in a plane diagram, \cite{Kauffman}, \cite{K15-1}, \cite{K15-2}, \cite{K17}, \cite{Satoh}, \cite{Turaev}. 
A {\it stabilization} of a surface-link $F$ is a connected sum 
$\bar F= F\# _{k=1}^s T_k$ of $F$ and a system of trivial torus-knots 
$T_k\, (k=1,2,\dots,s)$. 
By granting $s=0$, a surface-link $F$ itself is regarded as a stabilization 
of $F$. 
The trivial torus-knot system $T$ is called the {\it stabilizer} with {\it stabilizer components} $T_k\, (k=1,2,\dots,s)$ on the stabilization $\bar F$ of $F$. 
A {\it stable-ribbon} surface-link is a surface-link $F$ such that a stabilization 
$\bar F$ of $F$ is a ribbon surface-link. 
Every surface-link $F$ is equivalent to a stabilization of  a surface-link $F_*$ with minimal total genus.  This surface-link $F_*$ is called a {\it handle-irreducible summand} of $F$. 
The following result called {\it Stable-Ribbon Theorem} is our main theorem. 

\phantom{x}

\noindent{\bf Theorem~1.1.} 
Any handle-irreducible summand $F_*$ of every stable-ribbon surface-link $F$ is a ribbon 
surface-link which is determined uniquely from $F$ up to equivalences and stabilizations.

\phantom{x}

Any stabilization of a ribbon surface-link is a ribbon surface-link. So,  
 the following corollary is obtained from Theorem~1.1.

\phantom{x}

\noindent{\bf Corollary~1.2.} 
Every stable-ribbon surface-link is a ribbon surface-link. 

\phantom{x}

A {\it stably trivial} surface-link is a surface-link $F$ such that a stabilization 
$\bar F$ of $F$ is a trivial surface-link. 
Since a trivial surface-link is a ribbon surface-link, Theorem~1.1 also implies the following corollary, which is used to prove  smooth unknotting conjecture for a surface-link, \cite{K3}. 
This result leads to   4D smooth  and then classical  Poincar{\'e} conjectures,    
\cite{B}, \cite{K4}, \cite{K5-1},\cite{Per},  \cite{P1}, \cite{P2}.

\phantom{x}

\noindent{\bf Corollary~1.3.} 
Any handle-irreducible summand of every stably trivial surface-link is a trivial $S^2$-link, 
so that every stably trivial surface-link is a trivial surface-link.

\phantom{x} 

The plan for the proof of Theorem~1.1 is to show the following two lemmas by 
using the previous techniques,  \cite{K21}, \cite{K23}. 

\phantom{x}

\noindent{\bf Lemma~I.} Any handle-irreducible summand of any surface-link is unique up to equivalences and stabilizations.

\phantom{x}

\noindent{\bf Lemma~II.} Any stable-ribbon surface-link is a ribbon surface-link.

\phantom{x}

The proof of Theorem~1.1 is completed by these lemmas as follows:

\phantom{x}

\noindent{\bf Proof of Theorem~1.1 assuming Lemmas I, II.} By Lemma~II,  
any handle-irreducible summand of 
every stable-ribbon surface-link is a ribbon surface-link, which is unique up to equivalences 
and stabilizations by Lemma~I. This completes the proof of Theorem~1.1.

\phantom{x}

An idea of the proof of Lemma~I is to generalize the uniqueness result of an O2-handle pair 
on a surface-link earlier established to the case where the restriction on the attaching part is relaxed (see Theorem~2.2). 
An idea of the proof of Lemma~II is to consider a semi-unknotted multi-punctured handlebody system, simply called a {\it SUPH system}, of a ribbon surface-link. 
Two applications of Theorem~1.1 are made. One observation on Theorem~1.1 is the following theorem.

\phantom{x}

\noindent{\bf Theorem~1.4.} 
A connected sum $F=F_1\# F_2$ of surface-links $F_i\, (i=1,2)$ in $S^4$
is a ribbon surface-link if and only if both the surface-links 
$F_i\, (i=1,2)$ are  ribbon surface-links. 

\phantom{x}

This theorem contrasts with the behavior of  classical ribbon knot, 
because every classical knot is a connected summand of a connected sum ribbon knot. 
In fact, for every knot $k$ and the inversed mirror image $-k^*$ of $k$ in the 3-sphere $S^3$, 
the connected sum $k\#(-k^*)$ is a ribbon knot in $S^3$,   \cite{Fox1}, \cite{Fox2}, 
\cite{FM}, \cite{KSS-1}. 
A natural presentation of   $k\#(-k^*)$  is seen in a chord diagram  of  the spun 
$S^2$-knot of $k$ as a ribbon $S^2$-knot,  \cite{K15-1}. 

A surface-knot $F'$ in $S^4$ is obtained from a surface-link 
$F$ of $r$ components in $S^4$ by {\it fusion} if $F'$ is obtained from $F$ by surgery along $r-1$ disjointedly embedded 1-handles on $F$ in $S^4$. 
In earlier preprint of this paper, it is claimed that every $S^2$-link $L$ in $S^4$ consisting of trivial components in $S^4$ is a ribbon $S^2$-link. However, the proof contains an error. In fact, 
there was a non-ribbon $S^2$-link $L$ consisting of two trivial components such that 
the $S^2$-link obtained from $L$ by a fusion is  a non-ribbon $S^2$-knot, \cite{Ogasa}. 
As a revised content, the following theorem is shown, giving a characterization of when a surface-link consisting of ribbon surface-knot components is a ribbon surface-link.
Since a trivial surface-knot is a ribbon surface-knot, this theorem  also gives a characterization of when a surface-link consisting of trivial components is a ribbon surface-link.

\phantom{x}

\noindent{\bf Theorem~1.5.} The following statements (1)-(3) on a surface-link $F$ consisting of ribbon surface-knot components in $S^4$ are mutually equivalent.

\medskip

\noindent{(1)} $F$ is a ribbon surface-link.

\medskip

\noindent{(2)} The surface-knot obtained from  $F$ by every fusion  is a  ribbon surface-knot.

\medskip

\noindent{(3)} The surface-knot  obtained from  $F$ by a fusion  is a ribbon surface-knot.

\phantom{x}

Fundamental group is the most powerful invariant of a surface-link, \cite{CF}, \cite{Spanier}. 
Note that there are many $S^2$-links $L$ in $S^4$ whose fundamental groups 
$\pi_1(S^4\setminus L,x_0)$ have non-trivial torsion elements, \cite{Z}. 
These $S^2$-links are non-ribbon $S^2$-links because the fundamental group of 
a ribbon $S^2$-link is torsion-free, \cite {K5-2}. 
There are lots of classical non-ribbon links consisting of trivial components producing 
a ribbon-knot  by a fusion such as Hopf link, the split sum of Whitehead link and its mirror image, and the split sum of Borromean rings and its mirror image, etc., \cite{Fox1}, \cite{R}. Thus, this theorem also contrasts with the behavior of a classical ribbon link.

The proofs of Lemmas~I and II are given in Sections~2 and 3, respectively. 
In Section~4, the proofs of Theorems~1.4 and 1.5 are given.

\phantom{x}

\noindent{\bf 2. Proof of Lemma~I}

A 2-{\it handle} on a surface-link $F$ in ${\mathbf R}^4$ is a  2-handle 
$D\times I$ on $F$ with $D$ a core disk  embedded in   ${\mathbf R}^4$  such that 
$D\times I\cap F =\partial D\times I$, where 
$I$ denotes a closed interval containing $0$ and $D\times 0$ is identified with $D$. 
Two 2-handles $D\times I$ and $E\times I$ on $F$ are {\it equivalent} if there is an equivalence $f:{\mathbf R}^4\to {\mathbf R}^4$ from $F$ to itself 
such that the restriction $f|_F: F\to F$ is the identity map 
and $f(D\times I)=E\times I$. 

An {\it orthogonal 2-handle pair} (or simply, an {\it O2-handle pair}) on $F$ is a pair 
$(D\times I, D' \times I)$ of 2-handles 
$D\times I$, $D' \times I$ on $F$ such that 
\[D\times I\cap D' \times I = 
\partial D\times I\cap \partial D' \times I\]
and 
$\partial D\times I$ and $\partial D' \times I$ {\it meet orthogonally} on $F$, 
that is, the boundary circles $\partial D$ and $\partial D'$ meet transversely at one point $p$ so that the intersection $\partial D\times I\cap \partial D' \times I$ 
is homeomorphic to the square $Q=p \times I\times I$.
Let $(D\times I, D'\times I)$ be an O2-handle pair on a surface-link $F$. 
Let $F(D\times I)$ and $F(D'\times I)$ be the surface-links obtained from $F$ by the surgeries along $D\times I$ and $D'\times I$, respectively. 
Let $F(D\times I, D'\times I)$ be the surface-link which is 
the union $\delta\cup F^c_{\delta}$ of the plumbed disk 
\[\delta=\delta_{D\times I,D'\times I}
=D\times \partial I\cup Q\cup D'\times \partial I\]
and the surface 
\[F^c_{\delta}= \mbox{cl}(F\setminus (\partial D\times I\cup \partial D'\times I)).\] 
A once-punctured torus $T^o$ in a 3-ball $B$ 
is {\it trivial} if $T^o$ is smoothly and properly embedded in $B$ which splits 
$B$ into two solid tori. 
A {\it bump} of a surface-link $F$ is a 3-ball $B$ in ${\mathbf R}^4$ with 
$F\cap B=T^o$ a trivial once-punctured torus in $B$. 
Let $F(B)$ be a surface-link $F^c_B\cup \delta_B$ which is the union of  the surface 
$F^c_B=\mbox{cl}(F\setminus T^o )$ and a disk 
$\delta_B$ in the 2-sphere $\partial B$ with $\partial \delta_B =\partial T^o$.  
A {\it cellular move} of a compact  (possibly, bounded) surface $P$ in ${\mathbf R}^4$ 
is a compact surface $\tilde P$ such that the intersection $P^o=P\cap \tilde P$  is a  
once-punctured compact surface of  $P$ and $\tilde P$ with 
$d =\mbox{cl}(P \setminus P^o)$ and $\tilde d=(\tilde P \setminus P^o)$ disks in the interiors of $P$ and $\tilde P$, respectively such that  the union 
$d\cup \tilde d$ is a 2-sphere bounding a 3-ball smoothly embedded in ${\mathbf R}^4$ 
and not meeting the interior of $P^o$. 
Note that $F(B)$ is uniquely determined up to cellular moves on  the disk $\delta_B$ keeping $F^c_B$ fixed. For an O2-handle pair $(D\times I, D'\times I)$ on a surface-link $F$,
let $\Delta=D\times I\cup D'\times I$ is a 3-ball in ${\mathbf R}^4$ called 
the {\it 2-handle union}. 
Consider the 3-ball $\Delta$ as a Seifert hypersurface of 
the trivial $S^2$-knot $K=\partial \Delta$ in ${\mathbf R}^4$ to construct a 3-ball 
$B_{\Delta}$ obtained from $\Delta$ by adding an outer boundary collar. 
This 3-ball $B_{\Delta}$ is a bump of $F$, which we call the 
{\it associated bump} of the O2-handle pair $(D\times I, D'\times I)$. 
When the union of the 3-ball $\Delta$ and a boundary collar of $F^c_{\delta}$ are deformed into 
the 3-space ${\mathbf R}^3\subset {\mathbf R}^4$, this associated bump $B_{\Delta}$ is also considered as a regular neighborhood of $\Delta$ in ${\mathbf R}^3$. 
It is observed that 
an O2-handle unordered pair $(D\times I, D'\times I )$ on a surface-link $F$ 
is  constructed uniquely from any given bump $B$ of $F$ in ${\mathbf R}^4$ with 
$F(D\times I, D'\times I)\cong F(B)$, \cite{K21}. 
Further, for any O2-handle pair $(D\times I, D'\times I)$ on any surface-link $F$ and the associated bump $B$, there are identifications 
\[F(B) = F(D\times I, D'\times I) = F(D\times I) = F(D'\times I)\]
by equivalences which are attained by cellular moves on the disk $\delta=\delta_{D\times I, D'\times I}$ keeping $F^c_{\delta}$ fixed, \cite{K21}.  
A once-punctured torus $T^o$ in a 4-ball $A$ 
is {\it trivial} if $T^o$ is smoothly and properly embedded in $A$ and there is 
a solid torus $V$ in $A$ 
with $\partial V= T^o\cup \delta_A$ for a disk $\delta_A$ in the 3-sphere $\partial A$. 
A {\it 4D bump} of a surface-link $F$ is a 4-ball $A$ in ${\mathbf R}^4$ with 
$F\cap A=T^o$ a trivial once-punctured torus in $A$. 
A 4D bump $A$ is obtained from a bump $B$ of a surface-link $F$ by taking 
a bi-collar $c(B\times[-1,1])$ of $B$ in ${\mathbf R}^4$ with $c(B\times 0)=B$. 
The following lemma is proved by using a 4D bump $A$.

\phantom{x}

\noindent{\bf Lemma~2.1.} Let $(D\times I, D'\times I)$ be any O2-handle pair on any 
surface-link $F$ in ${\mathbf R}^4$, and $T$ a trivial torus-knot in ${\mathbf R}^4$ 
with any given spin loop basis $(e,e')$. 
Then there is an equivalence $f: {\mathbf R}^4\to {\mathbf R}^4$ from 
the surface-link $F$ to a connected sum $F(D\times I, D'\times I)\# T$ 
keeping $F^c_{\delta}$ fixed such that $f(\partial D)=e$ and $f(\partial D')= e'$.

\phantom{x}

\noindent{\bf Proof of Lemma~2.1.}
Let $A$ be a 4D bump associated with the O2-handle pair 
$(D\times I, D'\times I)$ on $F$. 
Let $\delta_A$  be a disk in the 3-sphere  $\partial A$ such that there is a solid torus $V$ in 
$A$  whose boundary is the union of the trivial once-punctured torus $P=F\cap A$ and the disk $\delta_A$. 
This solid torus $V$ induces an equivalence  
$f':({\mathbf R}^4, F)\to ({\mathbf R}^4, F(D\times I, D'\times I)\# T)$ 
sending $P$  to 
the connected summand $T^o$ of a connected sum $F(D\times I, D'\times I)\# T$ in $A$. 
Let  $(\tilde e, \tilde e')$ be  the spin loop basis of $T^o$ which is the image of 
the spin loop pair $(\partial D, \partial D')$ on $F$ under $f'$. 
There is an orientation-preserving 
diffeomorphism $g:{\mathbf R}^4 \to {\mathbf R}^4$ 
with $g|_{\mbox{cl}({\mathbf R}^4\setminus A)}=1$ such that 
$g(\tilde e, \tilde e')=(e,e')$ by the previous techniques, \cite{Hiro}, \cite{K21}. 
The composition $f=g f'$ is a desired equivalence. 
This completes the proof of Lemma~2.1.

\phantom{x}

A surface-link $F$ has {\it only unique O2-handle pair in the rigid sense} if 
for any O2-handle pairs $(D\times I, D'\times I)$ and $(E\times I, E'\times I)$ 
on $F$ with $(\partial D)\times I =(\partial E)\times I$ and 
$(\partial D')\times I =(\partial E')\times I$, there is an equivalence 
$f: {\mathbf R}^4 \to {\mathbf R}^4$ from $F$ to itself keeping $F^c_{\delta}$ fixed 
such that $f(D\times I)= E\times I$ and $f(D'\times I)= E'\times I$. 
It is proved that every surface-link $F$ has 
only unique O2-handle pair in the rigid sense,   \cite{K21}, \cite{K23}. 
A surface-link $F$ has {\it only unique O2-handle pair in the soft sense} if 
for any O2-handle pairs $(D\times I, D'\times I)$ and $(E\times I, E'\times I)$ 
on $F$ attached to the same connected component, say $F_1$ of $F$, there is an equivalence 
$f: {\mathbf R}^4 \to {\mathbf R}^4$ from $F$ to itself keeping $F^{(1)}= F\setminus F_1$ 
fixed such that 
$f(D\times I)= E\times I$ and $f(D'\times I)= E'\times I$. 
A surface-link not admitting any O2-handle pair is understood as 
a surface-link with only unique O2-handle pair in both the rigid and soft senses.
The following uniqueness of an O2-handle pair in the soft sense is essentially a consequence of the uniqueness of an O2-handle pair in the rigid sense.

\phantom{x}

\noindent{\bf Theorem~2.2 (Uniqueness of an O2-handle pair in the soft sense).} Every  
surface-link has only unique O2-handle pair in the soft sense. 

\phantom{x}

\noindent{\bf Proof of Theorem~2.2.} 
Let $(D\times I, D'\times I)$ and 
$(E\times I, E'\times I)$ be any two O2-handle pairs on a surface-link $F$ 
attached to the same connected component $F_1$ of $F$. 
Let $F^{(1)}=F\setminus F_1$. 
By Lemma~2.1,  there is an equivalences $f: {\mathbf R}^4 \to {\mathbf R}^4$ from $F$  to  the connected sum  
\[F(D\times I, D'\times I)\# T=F^{(1)}\cup \tilde F_1\# T\]
keeping  $F^{(1)}$ fixed and sending $F_1$ to $\tilde F_1\# T$, 
where $\tilde F_1=F_1(D\times I, D'\times I)$  and $T$ is a  trivial torus-knot in 
${\mathbf R}^4$.
Similarly, there is an equivalences $f': {\mathbf R}^4 \to {\mathbf R}^4$ from $F$  to  the connected sum  
\[F(E\times I, E'\times I)\# T'=F^{(1)}\cup \tilde F'_1\# T'\]
keeping $F^{(1)}$ fixed and sending $F_1$ to $\tilde F'_1\# T'$, 
where $\tilde F'_1=F_1(E\times I, E'\times I)$  and $T'$ is a  trivial torus-knot in 
${\mathbf R}^4$.
The diffeomorphism $g=f'f^{-1}: {\mathbf R}^4 \to {\mathbf R}^4$ is an equivalence from 
$F^{(1)}\cup \tilde F_1\# T$ to $F^{(1)}\cup \tilde F'_1\# T'$ keeping $F^{(1)}$ fixed. 
The connected sum
$\tilde F_1\# T$ is obtained from the split union $\tilde F_1+ T$ in ${\mathbf R}^4$ by surgery along an embedded 1-handle $h$ connecting a disk $d_1\subset \tilde F_1$ and a disk $d\subset T$, and the connected sum $\tilde F'_1\#T'$ is obtained from the split union 
$\tilde F'_1+T'$ in ${\mathbf R}^4$ by surgery along an embedded 1-handle $h'$ connecting a disk $d'_1\subset \tilde F'_1$ and a disk $d'\subset \tilde T$. Then there is a 4-ball $A$ in 
${\mathbf R}^4$ such that $T^o =A \cap (F^{(1)}\cup \tilde F_1\# T)$  is a trivial once-punctured torus of $T$ in $A$ with $d_1$ a disk bounded by the trivial knot 
$\partial T^o$ in the 3-sphere $\partial A$. Similarly, there is a 4-ball $A'$ in ${\mathbf R}^4$ such that $(T')^o =A' \cap (F^{(1)}\cup\tilde F'_1\# T')$ is a trivial once-punctured torus of 
$T'$ in $A'$ with $d'_1$ a disk bounded by the trivial knot $\partial (T')^o$ in the 3-sphere 
$\partial A'$.
It may be assumed that $g(\partial d_1)= \partial d'_1$ by sliding the attaching loop 
$g(\partial d_1)$  in $g(\tilde F_1\# T)$ and/or  the attaching loop 
$\partial d'_1$  in $\tilde F'_1\# T'$.  
Then it is assumed that  $g(T^o) = (T')^o$ (For a special case that 
$g(T^0)=\mbox{cl}(\tilde F'_1\# T'\setminus (T')^o)$, there is a deformation from $g(A)$ into  $A'$  to obtain $g(T^0)=(T')^o$).
Further, by Lemma~2.1, it is assumed that
$g(f(\partial D),f(\partial D')) = (f'(\partial E),f'(\partial E'))$. Then 
\[(f'(\partial D),f'(\partial D'))=g(f(\partial D),f(\partial D'))=(f'(\partial E),f'(\partial E')).\]
Since every surface-link has only unique O2-handle pair in the rigid sense, there is an equivalence $g': {\mathbf R}^4 \to {\mathbf R}^4$
from $F^{(1)}\cup \tilde F'_1\# T'$ to itself keeping $F^{(1)}$ fixed such that 
$g'(f'(D)\times I, f'(D')\times I) =(f'(E)\times I, f'(E')\times I)$.
The composite equivalence $g^*=(f')^{-1}g'gf: {\mathbf R}^4 \to {\mathbf R}^4$
is an equivalence from $F$ to itself keeping $F^{(1)}$ fixed and sending 
$(D\times I,D'\times I)$ to $(E\times I, E'\times I)$. Thus, every surface-link $F$ has only unique O2-handle pair in the soft sense. This completes the proof of Theorem~2.2.

\phantom{x}

The following corollary is obtained from the proof of Theorem~2.2.

\phantom{x}

\noindent{\bf Corollary~2.3.} Let $F$ and $ F'$ be surface-links with ordered components 
$F_i\, (i=1,2,\dots,r)$ and $ F'_i\, (i=1,2,\dots,r)$, respectively. 
Assume that the stabilizations
$\bar F= F\#_i T, \bar F'= F'\#_i T$ of $F, F'$ with induced ordered components obtained by the connected sums $F_i\# T, F'_i\# T$ of the $i$th components 
$F_i, F'_i$  and a trivial torus-knot $T$, respectively are equivalent by a component-order-preserving equivalence ${\mathbf R}^4 \to {\mathbf R}^4$. 
Then $F$ is equivalent to $F'$ by a component-order-preserving equivalence 
${\mathbf R}^4 \to {\mathbf R}^4$.

\phantom{x}

\noindent{\bf Remark~2.4.}  For the case of  ribbon surface-links $F$ and $F'$, Corollary~2.3  
has a different proof, \cite{K15-1}, \cite{K15-2}, \cite{K15-3}. 

\phantom{x}

The proof of Lemma~I is done as follows.

\phantom{x}

\noindent{\bf Proof of Lemma~I.} 
A surface-link $F$ with $r$ ordered components is $k$th 
{\it handle-reducible} if $F$ is equivalent to 
a stabilization $F'\#_k n_k T$ of a surface-link $F'$ for an integer $n_k>0$, where $\#_k n_k T$ denotes the stabilizer components $n_k T$ attaching to the $k$th component of $F'$. 
Otherwise, the surface-link $F$ is  said to be $k${\it th   handle-irreducible}. 
Note that if a surface-link $G$ is equivalent to a $k$th handle-irreducible surface-link $F$ by  
component-order-preserving equivalence, then $G$ is also 
$k$th handle-irreducible. 
Let $F$ and $G$ be ribbon surface-links with components $F_i\, (i=1,2,\dots,r)$ and $G_i\, (i=1,2,\dots,r)$, respectively. 
Let $F_*$ and $G_*$ be handle-irreducible summands of $F$ and $G$, respectively. 
Assume that there is an equivalence $f$ from $F$ to $G$. Then it is shown that 
$F_*$ and $G_*$ are equivalent as follows. 
Changing the indexes if necessary, assume that $f$ sends $F_i$ to $G_i$ for every 
$i$. 

Let 
\begin{eqnarray*}
F&=&F_*\#_1 n_1 T\#_2 n_2 T\#_3 \dots \#_r n_r T,\\
G&=&G_*\#_1 n'_1 T\#_2 n'_2 T\#_3 \dots \#_r n'_r T.
\end{eqnarray*}
 If necessary, by taking the inverse equivalence $f^{-1}$ instead of $f$,   assume 
 that $n'_1\geq n_1$. 
If $n'_1>n_1$, then  
there is a component-order-preserving equivalence $f^{(1)}$ from 
the first-handle-irreducible surface-link 
\[F_{(1)}=F_*\#_2 n_2 T\#_3 \dots \#_r n_r T\]
to the first-handle-reducible surface-link 
\[G_*\#_1 (n'_1-n_1) T\#_2 n'_2 T\#_3 \dots \#_r n'_r T,\]
by Corollary~2.3, which contradicts the first handle-irreducibility. 
Thus, $n'_1=n_1$ and the first handle-irreducible surface-link $F_{(1)}$ 
is equivalent to the first-handle-irreducible ribbon surface-link 
\[G_{(1)}=G_*\#_2 n'_2 T\#_3 \dots\#_r n'_r T.\]
By continuing this process, it is shown that $F_*$ is equivalent to $G_*$. 
This completes the proof of Lemma~I.

\phantom{x}

\noindent{\bf 3. Proof of Lemma~II}

A {\it chorded loop system} is a pair $(o,\alpha)$ of a trivial link $o$ and an arc system 
$\alpha$ attaching to $o$ in the 3-space ${\mathbf R}^3$, where $o$ 
and $\alpha$ are called a 
{\it based loop system} and a {\it chord system}, respectively. 
A {\it chorded loop diagram} or simply a {\it chord diagram} is a diagram $C(o,\alpha)$ in the plane ${\mathbf R}^2$ of the spatial graph $o\cup\alpha$. 
Let $D^+$ be a proper disk system in the upper half-space ${\mathbf R}_+^4$ 
obtained from a disk system $d^+$ in ${\mathbf R}^3$ 
bounded by $o$ by pushing the interior into ${\mathbf R}_+^4$. 
Similarly, let $D^-$ be a proper disk system in the lower half-space ${\mathbf R}_-^4$ 
obtained from a disk system $d^-$ in ${\mathbf R}^3$ 
bounded by $o$ by pushing the interior into ${\mathbf R}_-^4$. 
Let $O$ be the union of $D^+$ and $D^-$ which is a trivial $n S^2$-link in the 
4-space ${\mathbf R}^4$, where $n$ is the number of components of $o$. 
The union $O\cup \alpha$ is called 
a {\it chorded sphere system} constructed from a  chorded loop system $(o,\alpha)$. 
The chorded sphere system 
$O\cup \alpha$ up to orientation-preserving diffeomorphisms of ${\mathbf R}^4$ is 
independent of choices of $d^+$ and $d^-$ and uniquely determined by the chorded loop 
system $(o,\alpha)$ by  the Horibe-Yanagawa lemma,  \cite{KSS-1}. 
Thus, every ribbon surface-link  $F$ is uniquely constructed from 
a chorded loop system $(o,\alpha)$ via the chorded sphere system 
$O\cup \alpha$ so that $F=F(o,\alpha)$ is obtained from $O$ by surgery along a 1-handle system $N(\alpha)$ on $O$ with core arc system $\alpha$,  
where note  that the surface-link $F$ up to equivalences is unaffected by choices of any 
1-handle system $N(\alpha)$ fixing $\alpha$, \cite{HoK}.  
The moves on a chorded loop system $(o,\alpha)$ giving the same ribbon surface-link up to equivalences are determined, \cite{K15-1}, \cite{K15-2}. 
A {\it semi-unknotted multi-punctured handlebody system}  or simply a {\it SUPH system} 
for a surface-link $F$ in ${\mathbf R}^4$ is a  multi-punctured handlebody system $V$ 
(smoothly embedded) in ${\mathbf R}^4$ such that 
$\partial V=F\cup O$ for a trivial $S^2$-link $O$ in ${\mathbf R}^4$. Note that the numbers of 
connected components in $F$ and $V$ are equal.
The following lemma makes a characterization of a ribbon surface-link, 
\cite{KSS-2},  \cite{Yana}.

\phantom{x}

\noindent{\bf Lemma 3.1.} A surface-link $F$ is a ribbon surface-link 
if and only if $F$ admits a SUPH system $V$ in ${\mathbf R}^4$. 

\phantom{x}

\noindent{\bf Proof of Lemma~3.1.} A SUPH system 
$V$ for a ribbon surface-link $F$ is constructed from a chorded sphere system 
$O\cup \alpha$ by taking the union of a thickening $O\times [0,1]$ of $O$ in ${\mathbf R}^4$ 
and the 1-handle system $N(\alpha)$ attaching only to $O\times 0$. Conversely, given a SUPH system $V$ in ${\mathbf R}^4$ with $\partial V=F\cup O$ for a trivial $S^2$-link $O$, then take a chord system $\alpha$ in $V$ attaching to $O$
so that the frontier of the regular neighborhood of $O\cup \alpha$ in 
$V$ is parallel to $F$ in $V$. The chorded sphere system 
$O\cup \alpha$ shows that $F$ is a ribbon surface-link.
This completes the proof of Lemma~3.1.

\phantom{x}

Let $F$ be a surface-link of components $F_i\, (i=1,2,\dots, r)$ in ${\mathbf R}^4$. 
Let $F\# T$ be the connected sum of $F$ and 
a trivial torus-knot $T$ in ${\mathbf R}^4$ 
consisting of the components $F_1\# T, F_i\,(i=2,3,\dots,r)$. 
Assume that $F\# T$ is a ribbon surface-link. 
By Lemma 3.1, let $V$ be a SUPH system for $F\# T$ in ${\mathbf R}^4$. 
Let $V_1$ be the component of $V$ for $F_1\# T$ and write $V_1=U\#_{\partial} W$, a disk sum for a multi-punctured 3-ball $U$ and a handlebody $W$. 
The following lemma is needed to prove Lemma~II. 
 
\phantom{x}

\noindent{\bf Lemma~3.2.} For a suitable spin loop basis $(\ell, \ell')$ for $T^o$, there is 
a spin simple loop $ \tilde \ell'$ in the ribbon-surface-link $F_1\# T$ with intersection number 
$\mbox{Int}(\ell,\tilde \ell')\ne 0$ in $F_1\# T$ such that the loop 
$\tilde \ell'$ bounds a disk $D'$ in the handlebody $W$.

\phantom{x}

\noindent{\bf Proof of Lemma~3.2. } 
Consider a disk sum decomposition of the handlebody $W$ 
into solid tori $S^1\times D^2_j\, (j=1,2,\dots, g)$ pasting
along mutually disjoint disks. 
Let $(\ell_j,m_j)$ be a longitude-meridian pair of 
the solid torus $S^1\times D^2_j$ for all $j$. 
The loop basis 
$(\ell_j,m_j)$ for $S^1\times D^2_j$ is chosen to be a spin loop basis in ${\mathbf R}^4$ 
for all $j$,  \cite{HiK}.   By a choice of a spin loop basis $(\ell, \ell')$ for $T^o$, 
the loop $\ell$ meets a meridian loop $m_j$ with a non-zero intersection number 
in $\partial W$. The loop $m_j$ is taken to be a loop $\tilde \ell'$ in $F_1\#T$ 
bounding a disk $D'$ in $W$ with intersection number $\mbox{Int}(\ell,\tilde \ell')\ne 0$ since 
$m_j$ bounds a meridian disk $1\times D^2_j$ of the solid torus $S^1\times D^2_j\subset W$. 
This completes the proof of Lemma~3.2.

\phantom{x}

The following lemma is obtained by using Lemma~3.2.

\phantom{x}

\noindent{\bf Lemma~3.3.} 
There is a stabilization $\bar F$ of the ribbon surface-link $F\# T$ in ${\mathbf R}^4$ 
consisting of the components $\bar F_1, F_i\,(i=2,3,\dots,r)$ 
where $\bar F_1$ is the connected sum of $F_1\# T$ and trivial torus-knots $T_i\,(i=1,2,\dots, m)$ for some $m\geq 0$ such that the surface-link $\bar F$ has the following conditions 
(i) and (ii). 

\medskip

\noindent{(i)} There is an O2-handle pair $(D\times I, D'\times I)$ on 
$\bar F$ attached to $\bar F_1$ such that the surface-link 
$\bar F(D'\times I)$ is a ribbon surface-link with trivial 1-handles $h'_i\, (i=1,2,\dots, m)$ attached. 

\medskip

\noindent{(ii)} There is an O2-handle pair $(E\times I, E'\times I)$ on 
$\bar F$ attached to $\bar F_1$ such that the surface-link $\bar F(E'\times I)$ 
is $F$ with trivial 1-handles $h''_i\, (i=1,2,\dots, m)$ attached.

\phantom{x}

\noindent{\bf Proof of Lemma~3.3.}
Let $p_i\, (i=0,1,\dots, m)$ be the intersection points of 
transversely meeting simple loops $\ell$ and $\tilde\ell'$ in $F_1\# T$ given by Lemma~3.2.
For every $i>0$, let $\alpha_i$ be an arc neighborhood of $p_i$ in $\ell$, 
and $h_i$ a 1-handle on $F\# T$ with a core arc $\hat\alpha_i$ obtained 
by pushing the interior of $\alpha_i$ outside the SUPH system $V$.  
Let $\bar F= F\#T \#_{i=1}^m T_i$ be a stabilization of $F\# T$ with the component 
$\bar F_1= F_1\#T \#_{i=1}^m T_i$ obtained from $F_1\# T$ 
by surgery along the disjoint trivial 1-handle system $h_i\, (i=1,2,\dots, m)$. 
Let  $\alpha_i^+=\alpha_i\cup (h_i\cap \ell) $ be the arc in  
$\ell$ extending $\alpha_i$. 
Let $\tilde \alpha_i$ be a proper arc in the annulus 
$\mbox{cl}(\partial h_i\setminus h_i\cap F\#T)$ which is  parallel  
to the core arc $\hat \alpha_i$ of  $h_j$ in  $h_i$ with 
$\partial \tilde \alpha_i=\partial\alpha_i^+$. 
Let $\tilde\ell$ be a simple spin loop in $\bar F$ obtained from $\ell$ 
by replacing $\alpha_i^+$ with $\tilde \alpha_i$ for every $i>0$, which  
meets $\tilde\ell'$ transversely in just one point. 
Let $W^+(D')$ be the handlebody obtained from the handlebody 
$W^+=W \cup_{i=1}^{m} h_i$ by removing a thickened disk $D'\times I$ 
of $D'$. The manifold $V^+(D')$ obtained from the SUPH system 
$V^+=V\cup_{i=1}^{m} h_i$ for the ribbon surface-link $\bar F$  
by replacing  $W^+$ with $W^+(D')$ is 
a SUPH system for a surface-link $\bar F'$ in ${\mathbf R}^4$ consisting of  $F_i\,(i=2,3,\dots,r)$ and a component $\bar F'_1$ with genus reduced by $1$ from $\bar F_1$. 
By Lemma~3.1, $\bar F'$ is a ribbon-surface-link in ${\mathbf R}^4$. 
The SUPH system $V^+$ for $\bar F$ is a disk sum of $V^+(D')$ and  a solid torus $W_1$ 
with  the disk $D'$ as a meridian disk  and  the loop $\tilde\ell$ as a longitude. 
Let $d_W=V^+(D')\cap W_1$ be the pasting disk between $V^+(D')$ and $W_1$, which is regarded as a 1-handle $h_W$ joining $V^+(D')$ and $W_1$. 
Let $(E\times I, E'\times I)$ be an O2-handle pair  on $F\# T$ in ${\mathbf R}^4$ attached to 
$T^o$ with $(\partial E, \partial E')=(\ell, \ell')$ in Lemma~3.2. 
Let $A$ be a 4D bump of the associated bump $B$ of $(E\times I, E'\times I)$. 
In the case of (i), since there is no need to worry about the intersection of $A$ with 
$E\cup E'$, the 4D ball $A$ is deformed so that $V^+\cap A=W_1\cup h_W$ 
by observing that $V\setminus V_1$ and $U$ are  disjoint from $A$ by construction of $A$ 
and by taking spine graphs of $W^+(D')$, $W_1$ and $h_W$. 
Then the loop $\tilde\ell$ bounds a disk $D$ in $A$ 
not meeting the interior of $W_1$ and $h_W$. 
This means that there is an O2-handle pair 
$(D\times I, D'\times I)$ on the surface-link $\bar F$ such that 
$\bar F(D'\times I)$ is a ribbon surface-link with 
trivial 1-handles $h'_i\, (i=1,2,\dots,m)$ attached, showing (i). 
For the case of (ii), note that the 1-handles $h_i\, (i=1,2,\dots,m)$ on $F\# T$ are  deformed  isotopically in $A$ into 1-handles $h''_i\, (i=1,2,\dots,m)$ on $F\# T$ disjoint from the disk pair 
$(E,E')$ because the core arcs of the 1-handles $h_i\, (i=1,2,\dots,m)$ are deformed to be 
disjoint from the disk pair $(E,E')$ in $A$. 
The surface-link $\bar F(E'\times I)$ which is equivalent to $\bar F(E\times I, E'\times I)$ 
is the surface-link $F$ with the trivial 1-handles 
$h''_i\, (i=1,2,\dots,m)$ attached, showing (ii). 
Thus, the proof of Lemma~3.3 is completed.

\phantom{x}

The following lemma is a combination of Lemma~3.3 and the uniqueness of an O2-handle pair in the soft sense (Theorem~2.2). 

\phantom{x}

\noindent{\bf Lemma~3.4.} 
If a connected sum $F\# T$ of a surface-link $F$ and a trivial torus-knot $T$ 
in ${\mathbf R}^4$ is a ribbon surface-link, then $F$ is a ribbon surface-link.

\phantom{x}

\noindent{\bf Proof of Lemma~3.4.}
Let $F\# T=F_1\# T\cup F_2\cup\dots\cup F_r$ be a ribbon surface-link 
for a trivial torus-knot $T$. 
By Lemma~3.3 (i), the surface-link $F''=\bar F(D\times I, D'\times I)$ equivalent to  
$\bar F(D\times I)$ is a ribbon surface-link and further the surface-link 
$F^*$ obtained from $F''$ by the surgery on O2-handle pairs of all the trivial 
1-handles $h'_i\, (i=1,2,\dots, m)$ is  a ribbon surface-link. 
By Lemma~3.3 (ii), the surface-link $\bar F(E\times I, E'\times I)$ equivalent to 
$\bar F(E\times I)$
is the surface-link $F$ with the 1-handles $h''_i\, (i=1,2,\dots, m)$ trivially attached. 
By an inductive use of Theorem~2.2 (Uniqueness of an O2-handle pair in the soft sense), 
the surface-link $F$ is equivalent to the ribbon surface-link $F^*$. 
Thus, $F$ is a ribbon surface-link and the proof of Lemma~3.4 is completed. 

.
\phantom{x}

Lemma~II is a direct consequence of Lemma~3.4 as follows. 

\phantom{x}

\noindent{\bf Proof of Lemma~II.} If a stabilization 
$\bar F$ of a surface-link $F$ is a 
ribbon surface-link, then $F$ is a ribbon surface-link by an inductive use 
of Lemma~3.4. This completes the proof of Lemma~II.

\phantom{x}

\noindent{\bf 4. Proofs of Theorems~1.4 and 1.5}

The proof of Theorem~1.4 is done as follows.

\phantom{x}

\noindent{\bf Proof of Theorem~1.4.} The \lq if\rq\, part of Theorem~1.4 
is seen from the definition of a ribbon surface-link. 
The proof of the \lq only if\rq\, part of Theorem~1.4 uses 
the fact that every surface-link is made a trivial surface-knot by  surgery along a finite number of possibly non-trivial 1-handles, \cite{HoK}. 
The connected summand $F_2$ of $F_1\# F_2$ is made a trivial surface-knot by  surgery 
along 1-handles within the 4-ball defining the connected summand $F_2$, so that 
the surface-link $F$ changes into a new ribbon surface-link and hence 
$F_1$ is a stable-ribbon surface-link. 
By Corollary~1.2, $F_1$ is a ribbon surface-link. 
By interchanging the roles of $F_1$ and $F_2$, the connected summand $F_2$ 
is also a ribbon surface-link. 
This completes the proof of Theorem~1.4.

\phantom{x}

The following lemma is not used in the present version of Theorem~1.5, but this lemma 
remains here since it is an interesting property. 

\phantom{x}

\noindent{\bf Lemma~4.1.}  Let $K$ be an $S^2$-knot  in $S^4$ obtained 
from a trivial surface-knot $F$ of genus $n$ in $S^4$  by surgery along disjoint 2-handles 
$D_i\times I\, (i=1,2,\dots,n)$. Then there is a disjoint O2-handle pair system
$(D_i\times I, D'_i\times I)  \, (i=1,2,\dots,n)$ on $F$ in $S^4$ if and only if $K$ is a 
trivial $S^2$-knot in $S^4$.  

\phantom{x}

\noindent{\bf Proof of Lemma~4.1.} If there is a  disjoint O2-handle pair system   
$(D_i\times I, D'_i\times I)  \, (i=1,2,\dots,n)$ on $F$ in $S^4$, then $K$ is a trivial $S^2$-knot 
by Corollary~1.3. 
Note that the 2-handle system   $D_i\times I\, (i=1,2,\dots,n)$ on $F$ is a 1-handle system 
on $K$. If $K$ is a trivial $S^2$-knot, then the 1-handle system $D_i\times I\, (i=1,2,\dots,n)$ 
on the trivial $S^2$-knot $K$ is always a trivial 1-handle system on $K$, \cite{HoK}.
Hence there is  a disjoint O2-handle pair system 
$(D_i\times I, D'_i\times I)  \, (i=1,2,\dots,n)$ on $F$ in $S^4$. 
This completes the proof of Lemma~4.1.

\phantom{x}

The proof of Theorem~1.5 is done as follows.

\phantom{x}

\noindent{\bf Proof of Theorem~1.5.} 
The assertions (1) $\to$ (2) and (2) $\to$ (3) are obvious by definitions. 
The assertion (3) $\to$ (1) is shown as follows. 
When $F$ is a ribbon surface-knot, there is nothing to prove. 
By inductive assumption, the result for $F$ of $r-1$ components is assumed. 
Let $F$ be a surface-link of $r$ ribbon surface-knot components $F_i,\,(i=1,2,\dots,r)$. 
Since a fusion of $F$ makes a ribbon surface-knot by assumption, 
let $K$ be a ribbon surface-knot obtained form 
$F$ by a fusion along a disjoint 1-handle system $h$ on $F$ with only one 1-handle $h_r$ connecting to the component $F_r$. 
Let $K'$ be a surface-knot obtained from the sublink $F\setminus F_r$ by the fusion 
along the 1-handle system $h\setminus h_r$. 
Since $F_r$ is a ribbon surface-knot, let $(O_r,\alpha_r)$ be a chorded sphere-system for 
the ribbon surface-knot $F_r$. Since $[K']=0$ in $H_2(S^4\setminus O_r;Z)=0$, there is a compact connected oriented 3-manifold $V'$ in $S^4$ with $\partial V'=K'$ and $V'\cap O_r=\emptyset$. 
Since the arc system $\alpha_r$ transversely meets with the interior of $V'$, a multi-punctured 
manifold $(V')^{(0)}$ of $V'$ does not meet the chorded sphere system $O_r\cup\alpha_r$, so that $W_r\cap (V')^{(0)}=\emptyset$ for a SUPH system $W_r$ for $F_r$ with 
$\partial W_r =F_r\cup O_r$, constructed from $O_r\cup\alpha_r$. 
Let $g$ be a disjoint 1-handle system on $F'$ embedded in $(V')^{(0)}$ such that the closed complement $H^{(0)}=\mbox{cl}((V')^{(0)}\setminus g)$
is a multi-punctured handlebody of genus, say $n$. 
Let $\partial H^{(0)}=(K')^+ \cup O^H$ where $(K')^+$ is the surface-knot obtained from 
$K'$ by surgery along the 1-handle system $g$. 
The union $W= H^{(0)}\cup W_r$ is a SUPH system for the surface-link 
$(K')^+\cup F_r$ in $S^4$, so that $(K')^+\cup F_r$ is a ribbon surface-link in $S^4$ by Lemma~3.1. 
By replacing $W$ with a multi-puncture manifold of $W$, the union 
$W^+=W\cup h_r$ is a SUPH system for the ribbon surface-knot $K^+$ in $S^4$ obtained 
from the surface-link $(K')^+\cup F_r$ by fusion along $h_r$. Note that 
the ribbon surface-knot $K^+$ is also obtained from the ribbon surface-knot $K$ 
by surgery along $g$. 
Let $W_K$ be a SUPH system for the ribbon surface-knot $K$ in $S^4$. 
By replacing $W_K$ with a multi-punctured manifold of $W_K$, the union 
$W_K^+=W_K\cup g$ is a SUPH system for the ribbon surface-knot $K^+$ in $S^4$. 
Equivalent ribbon surface-links are faithfully equivalent and they are
moved into each other by the moves M0, M1, M2, \cite{K15-3}. This means that there is an
orientation-preserving diffeomorphism $f$ of $S^4$ sending 
a multi-punctured manifold 
of the SUPH system $(W^+)^{(0)}$ of $W^+$ to a multi-punctured manifold 
$(W_K^+)^{(0)}$ of the SUPH system $W_K^+$ and keeping $K^+$ set-wise fixed. 
Let $D(h_r)$ and $D(g)$ be a proper disk and a proper disk system in $h_r$ and $g$ 
parallel to the attaching disks with one disk for every 1-handle, respectively. 
The proper disk $f(D(h_r))$ and the proper disk system $D(g)$ may meet transversely 
with simple loops in the interior of the punctured handlebody $(W_K^+)^{(0)}$. 
Let $f(D(h_r))'$ be a proper disk disjoint from $D(g)$ with 
$\partial f(D(h_r))'= \partial f(D(h_r))$ in $(W_K^+)^{(0)}$ 
obtained from $f(D(h_r))$ by using a cutting technique 
along an innermost loop of $f(D(h_r))\cap D(g)$ in $D(g)$ inductively. 
By splitting $(W^+)^{(0)}$ along the disk system $f(D(h_r))'\cap D(g)$, a
SUPH system $W(K')\cup W(F_r)$ for the surface-link $K'\cup F_r$ is obtained, 
meaning that the surface-link $K'\cup F_r$
is a ribbon surface-link. In particular, $K'$ is a ribbon surface-knot in $S^4$. 
By inductive assumption, $F\setminus F_r$ is a ribbon surface-link. 
Let $W(F\setminus F_r)$ be a SUPH system for the ribbon surface-link $F\setminus F_r$ 
in $S^4$. 
A multi-punctured manifold of the SUPH system $W(F\setminus F_r)\cup (h\setminus h_r)$ 
for the ribbon surface-knot $K'$ is sent to a
multi-punctured manifold of the SUPH system $W(K')$ for $K'$ 
by an orientation-preserving diffeomorphism $f'$ of $S^4$. 
After replacing $W(F\setminus F_r)\cup (h\setminus h_r)$ and $W(K')\cup W(F_r)$ 
with multi-punctured manifolds, respectively, 
the preimage $(f')^{-1}(W(K')\cup W(F_r))$ is the union of 
$W(F\setminus F_r)\cup (h\setminus h_r)$ and $(f')^{-1}(W(F_r))$, showing that 
$W(F\setminus F_r)\cup (f')^{-1}(W(F_r))$ is a SUPH system for the surface-link $F$. 
Thus, $F$ is a ribbon surface-link. 
By induction on $r$, (3) $\to$ (1) is obtained and the proof of Theorem~1.5 is completed.

\phantom{x}

\noindent{\bf 5. Conclusion}

Ribbonness of a stable-ribbon surface-link shown in Theorem~1.1 is applied to determine the Ribbonness of some classes of surface-links as it is shown in Theorem~1.4. 
Theorem~1.5 seems unrelated to Theorem~1.1, but 
the result that equivalent ribbon surface-links are faithfully equivalent, which was used in the proof of Theorem~1.5, is obtained by using a result that is a predecessor of Theorem~1.1, 
\cite{K15-3}. 
Due to the error mentioned before the statement of Theorem~1.5, 
Theorem~1 of \cite {K6} should be replaced by Theorem~1.5 for a surface-link consisting of 
trivial components. Theorem~2 of \cite {K6} is fine as it is, which says that 
if $\mathbf F$ contains at least two non-sphere components, then 
there exists a pair of a ribbon $\mathbf F$-link and a non-ribbon $\mathbf F$-link 
with isomorphic fundamental groups by a meridian-preserving isomorphism.
Ribbonness of a surface-link relates not only to smooth unknotting conjecture for a surface-link leading to classical and 4D smooth Poincar{\'e} conjectures, but also to J. H. C. Whitehead asphericity conjecture for aspherical 2-complex, \cite{Hw}, \cite{K5-2}, \cite{K5-4}, \cite{Ker}, \cite{MKS}, \cite{Wh}, and Kervaire conjecture on group weight, \cite{GR}, \cite{K5-3}, \cite{Kly}. In these arguments, it is a basic result that every $S^2$-link with free 
fundamental group is a ribbon $S^2$-link, \cite{K5-2}, \cite{K5-5}, \cite{K5-6}. 
In another direction, it may be an interesting problem to investigate a canonical relationship between a chorded loop diagram of a ribbon surface-knot and a knot diagram, \cite{K5-1}. In conclusion, ribbon surface-knot theory will be a tool for studies of low dimensional topology.

\phantom{x}

\noindent{\bf Acknowledgements.} 
The author would like to thank O. Chterental (a student of D. Bar-Natan) for letting him know 
the reference \cite{Ogasa} and comments, which made the revision of Theorem~1.5 possible.
This work was partly supported by JSPS KAKENHI Grant Number JP21H00978 and MEXT Promotion of Distinctive Joint Research Center Program JPMXP0723833165.

\phantom{x}

\end{document}